\documentclass[multphys,vecphys]{svmult}

\usepackage{amssymb}
\usepackage{amsmath}
\usepackage{latexsym}
\usepackage{amscd}

\def\I{{\mathbb I}}
\def\R{{\mathbb R}}
\def\K{{\mathbb K}}
\def\C{{\mathbb C}}

\def\N{{\mathbb N}}
\def\H{{\mathcal H}}
\def\F{{\mathcal F}}
\def\T{{\mathcal T}}

\def\newblock{ }
\def\ds{\displaystyle}
\begin{document}

\title*{Deformation quantization on a Hilbert space}

\author{Giuseppe Dito}

\institute{Institut de Math\'ematiques de Bourgogne\\ Universit\'e de Bourgogne\\
UMR CNRS 5584\\ B.P. 40780, 
21078 Dijon Cedex, France\\
\texttt{giuseppe.dito@u-bourgogne.fr} }

\maketitle

\abstract{
We study deformation quantization on an infinite-dimensional
Hilbert space $W$ endowed with its canonical Poisson structure. 
The  standard example of the Moyal star-product is made explicit
and it is shown that it is well defined on a subalgebra of $C^\infty(W)$.
A classification of inequivalent deformation quantizations of exponential type,
containing the Moyal and normal star-products, is also given.}

\section{Introduction}\label{sec:intro}

Deformation quantization provides an alternative formulation of Quantum
Mechanics by interpreting quantization as a deformation of the
commutative algebra a classical observables into a noncommutative
algebra \cite{bayen.et.al:1978a}. The quantum algebra is defined by
a formal associative star-product $\star_\hbar$ which
encodes the algebraic structure of the set of observables.

Deformation quantization has been applied with increasing generality
to several area of mathematics and physics. Most of these applications
deal with star-products on finite-dimensional manifolds. 
See \cite{dito.sternheimer:2002a} for a recent review.

It is natural to consider an extension of deformation quantization to
infinite-dimensional manifolds as it appears to be a good setting where
quantum field theory of nonlinear wave equations 
can be formulated (e.g. in the sense of I.~Segal \cite{segal:1974a}).
In the star-product approach, the first steps in that direction are given in \cite{dito:1990a,dito:1992a}.

Recently, deformation quantization has become popular among  field and string theorists.
A generalization of Moyal star-product to infinite-dimensional spaces
appears in several places in the literature. 
Let us just notice that the  Witten star-product \cite{witten:1986} appearing
in string field theory is {\em heuristically} equivalent to
an infinite-dimensional version of the Moyal star-product.
A brute force generalization of Moyal star-product to field theory
yields to some pathological and unpleasant features
such as anomalies and breakdown of associativity.
We think that it is worth writing down a mathematical study of the Moyal
product in infinite dimension even if it is not an adequate product
for field theory considerations.

In the finite-dimensional case, the existence of star-products on any (real) symplectic manifold  
has been established by DeWilde and Lecomte \cite{dewilde.lecomte:1983b}. The general existence and classification
problems for the deformation quantization of a Poisson manifold was solved by Kontsevich \cite{kontsevich:1997}.
However, the very first problem that one faces when going over infinite-dimensional spaces,
it to make sense of the star-product itself as a formal associative product.
It contrasts with the finite-dimensional case where the deformation is defined
on all of the smooth functions on the manifold. This is by far too demanding 
in the infinite-dimensional case even when the Poisson structure is well-defined 
on all of the smooth functions (e.g. on Banach or Fr\'echet spaces). 
One should specify first an Abelian algebra of admissible functions which then
can be deformed. For example, on $E=\mathcal{S}\times\mathcal{S}$, where $\mathcal{S}$ is the
Schwartz space on $\R^n$, endowed with its canonical Poisson structure,
one cannot expect to write down a star-product defined on all holomorphic
functions on $E$, but has to restrict it to some subalgebra. For example, in \cite{dito:1992a} it is shown,
that for such a simple star-product as the normal star-product, it is defined on the subalgebra
of holomorphic functions of $a$ and $\bar a$ (creation and annihilation `operators') having
semi-regular kernels. In \cite{duetsch.fredenhagen:2001b}, one can find a nice  analysis for the normal 
star-product and the conditions on the kernel have been translated in terms of  wave front set of the distributions.

After making precise what is a deformation quantization on a Hilbert space, 
we first present a study of Moyal product when the space-space 
is the direct sum of a Hilbert space with its dual.
We identify a subalgebra of smooth functions, specified by conditions of Hilbert-Schmidt 
type on their derivatives, on which the Moyal product makes sense.
We also define a family of star-products of exponential type, show that
they are not all equivalent to each other and give the classification
of their equivalence classes in terms of Hilbert-Schmidt operators.

\section{Star-products on a Hilbert space}

When  infinite-dimensional spaces are involved, further conditions
are needed to define a deformation quantization or a star-product.
The algebra of functions on which the Poisson bracket and the star-product are 
defined should be specified along with the class of {\em admissible} cochains
(especially when the issue of the equivalence of deformations is considered). 

\subsection{Notations}

Let $B$ be a Banach space over a field $\K$ ($\R$ or $\C$). The topological dual of $B$
shall be denoted by $B^*$.
The Banach space  of bounded $r$-linear  forms on~$B$ is denoted by $\mathcal{L}^r(B,\K)$
and $\mathcal{L}^r_{\rm sym}(B,\K)$ is the subspace of $\mathcal{L}^r(B,\K)$ consisting
of bounded symmetric $r$-linear forms on $B$.
We shall denote by $C^\infty(B,\K)$ the space of $\K$-valued functions on $B$ that 
are smooth in the Fr\'echet sense.  The Fr\'echet derivative of $F\in C^\infty(W,\K)$ is denoted by $DF$ and 
it is a smooth map from $B$ to $\mathcal{L}^1(B,\K)=B^*$, i.e., $DF\in C^\infty(B,B^*)$.
For $F\in C^\infty(B,\K)$, the higher derivative $D^{(r)}F$ belongs to $C^\infty(B,\mathcal{L}^r_{\rm sym}(B,\K))$
and we shall use the following notation $D^{(r)}F(b).(b_1,\ldots,b_r)$ for the $r^{\rm th}$-derivative
of $F$ evaluated at $b\in B$ in the direction of $(b_1,\ldots,b_r)\in B^r$.

Let $W$ be an infinite-dimensional separable Hilbert space over a field  $\K$. 
For notational reasons, as it will become clear later, it would be convenient 
for us to not identify $W^*$ with $W$.
For any orthonormal basis $\{e_i\}_{i\geq1}$ in $W$ and corresponding
dual basis $\{e^*_i\}_{i\geq1}$ in $W^*$,
we shall denote the partial derivative of $F\in C^\infty(W,\K)$ evaluated at $w$ 
in the direction of $e_i$ by $\partial_iF(w)\in \K$, i.e., $\partial_iF(w)=DF(w).e_i$. 
Since $F$ is differentiable in the  Fr\'echet sense, we have $DF(w) = \sum_{i\geq1} \partial_iF(w)e^*_i$
and thus, for any $w\in W$, that $\sum_{i\geq1} |\partial_iF(w)|^2<\infty$.

\subsection{Multidifferential operators}

Let us first make precise what we call a Poisson
structure on $W$. In the following, we will consider  a map $P$ that sends $W$ into a space
of (not necessarily bounded) bilinear forms on $W^*$ and a subalgebra $\F$ of $C^\infty(W,\K)$.
We define the subspace $\mathcal{D}^\F_w =\{DF(w)\ | \ F\in \F \}$ of  $W^*$.
\vskip2mm
\begin{definition}\label{def:poisson}
Let $W$ be a Hilbert space. Let $\F$ be an Abelian subalgebra (for the pointwise product)
of $C^\infty(W,\K)$. 
A Poisson bracket on $(W,\F)$ is a $\K$-bilinear map $\{\cdot,\cdot\}\colon \F\times \F \rightarrow \F$ 
such that:
\vskip2mm
\noindent i) there exists a map $P$ from $W$ to the  space of bilinear forms on $W^*$,
so that the domain of $P(w)$ contains $\mathcal{D}^\F_w\times \mathcal{D}^\F_w$ and
$\forall F,G\in \F$, $\{F,G\}(w) = P(w).(DF(w),DG(w))$ where $w\in W$.
\vskip2mm
\noindent ii) $(\F, \{\cdot,\cdot\})$ is a Poisson algebra, 
i.e., skew-symmetry, Leibniz rule, and Jacobi identity are satisfied.
\vskip2mm
\noindent The triple $(W,\F,\{\cdot,\cdot\})$ is called a Poisson space.
\end{definition}
\vskip2mm
Let us give an example where $P(w)$ is an unbounded bilinear form on $W^*$.

\begin{example}
Consider a real Hilbert space $W$ with orthonormal basis $\{e_i\}_{i\geq0}$.
We will realize the following subalgebra of the Witt algebra:
$$
[L_m,L_n] = (m-n) L_{m+n},\quad m,n\geq 0,
$$
by functions on $W$. For $w\in W$, let $\phi_i(w)=\langle e^*_i, w\rangle$, $i\geq0$ be the coordinate
functions. The algebra $\mathcal F$ generated by the family of functions $\{\phi_i\}_{i\geq0}$ 
is an Abelian  subalgebra of $C^\infty(W,\R)$ consisting of polynomial functions in a finite number
of variables. The following expression:
$$
\{F,G\}(w) = \sum_{m,n\geq0} (m-n) \phi_{m+n}(w) \partial_m F(w) \partial_n G(w),\quad F,G\in \F, w \in W,
$$
defines a Poisson bracket on $(W,\F)$. Indeed the right-hand side is a finite sum and is a function
in  $\F$, and we have:
$$
\{\phi_i,\phi_j\}=(i-j)\phi_{i+j},
$$
from which Jacobi identity follows. The special case when $j=0$ gives:
$$
\{\phi_i,\phi_0\}(w)=P(w).(D\phi_i(w), D\phi_0(w)) = P(w).(e^*_i,e^*_0)=  i \phi_{i}(w).
$$
By choosing an appropriate $w$ (e.g. $w=\sum_{i\geq1} i^{-3/4} e_i$), then 
$P(w).(e^*_i,e^*_0)=i \phi_{i}(w)$ can become as large
as desired by varying $i$. This shows that the bilinear form $P(w)$ cannot be bounded. 
\end{example}

The generalization of Def.~\ref{def:poisson} to multidifferential operators on $W$ is straightforward.
Again, given an Abelian subalgebra $\F$ of $C^\infty(W,\K)$, we define the following subspace of
$\mathcal{L}^r_{\rm sym}(W,\K)$:
$$
\mathcal{D}^\F_w(r) =\{D^{(r)}F(w)\ | \ F\in \F \}.
$$
\vskip2mm
\begin{definition}\label{def:op}
Let $W$ be a Hilbert space. Let $\F$ be an Abelian subalgebra, for the pointwise product,
 of $C^\infty(W,\K)$. 
Let $r\geq1$, an $r$-differential operator ${\bf A}$ on $(W,\F)$ is an 
$r$-linear map ${\bf A} \colon \F^r \rightarrow \F$ 
such that:
\vskip2mm
\noindent i) for $(n_1,\ldots,n_r)\in\N^r$, there exists a  map $a^{(n_1,\ldots,n_r)}$ from 
$W$ to a space of (not necessarily bounded) $r$-linear forms on 
$\mathcal{L}^{n_1}_{\rm sym}(W,\K)\times\cdots\times \mathcal{L}^{n_r}_{\rm sym}(W,\K)$, i.e.,
$$
a^{(n_1,\ldots,n_r)}(w)\colon \mathcal{D}_w^{(n_1,\ldots,n_r)}\subset
\mathcal{L}^{n_1}_{\rm sym}(W,\K)\times\cdots\times \mathcal{L}^{n_r}_{\rm sym}(W,\K)\rightarrow \K,
$$
so that the domain $\mathcal{D}_w^{(n_1,\ldots,n_r)}$ of $a^{(n_1,\ldots,n_r)}(w)$ contains 
$\mathcal{D}^\F_w(n_1)\times\cdots \times\mathcal{D}^\F_w(n_r)$
and $a^{(n_1,\ldots,n_r)}$ is $0$ except for finitely many $(n_1,\ldots,n_r)$;
\vskip2mm
\noindent ii) for any $F_1,\ldots,F_r\in \F$ and $w\in W$, we have
$$
{\bf A}(F_1,\ldots,F_r)(w)=\sum_{n_1,\ldots,n_r\geq 0} a^{(n_1,\cdots ,n_r)}(w).(D^{(n_1)}F_1(w),\cdots,D^{(n_r)}F_r(w)).
$$
\end{definition}
\vskip2mm
Notice that Poisson brackets as defined above are special cases of 
bidifferential operators in the sense of Def.~\ref{def:op} with $P=a^{(1,1)}$.

\subsection{Deformation quantization on $W$}

We now have all the ingredients to define what is meant by deformation quantization
of a Poisson space $(W,\F,\{\cdot,\cdot\})$ when $W$ is a Hilbert space.
\vskip2mm
\begin{definition}
Let $W$ be a Hilbert space and $(W,\F,\{\cdot,\cdot\})$ be a Poisson space. 
A star-product on $(W,\F,\{\cdot,\cdot\})$ is a $\K[[\hbar]]$-bilinear product
$\star_\hbar\colon \F[[\hbar]]\times\F[[\hbar]]\rightarrow\F[[\hbar]]$ given by
$F\star_\hbar G = \sum_{r\geq0} \hbar^r C_r(F,G)$ for $F,G\in\F$ and extended by $\K[[\hbar]]$-bilinearity
to $\F[[\hbar]]$, and satisfying for any
$F,G,H\in\F$:
\vskip2mm
\noindent i) 
$C_0(F,G)= FG$,
\vskip2mm
\noindent ii) 
$C_1(F,G) - C_1(G,F) = 2 \{F,G\}$,
\vskip2mm
\noindent iii) 
for $r\geq1$, $C_r\colon\F\times\F\rightarrow \F$ are bidifferential operators in the sense of Def.~\ref{def:op},
vanishing on constants,
\vskip2mm
\noindent iv) 
$F\star_\hbar(G\star_\hbar H)=(F\star_\hbar G)\star_\hbar H$.
\vskip2mm
\noindent The triple $(W,\F[[\hbar]],\star_\hbar)$ is called a deformation quantization
of the Poisson space $(W,\F,\{\cdot,\cdot\})$.
\end{definition}
\vskip2mm
We also have a notion of equivalence of deformations adapted to our context:
\vskip2mm
\begin{definition}
Two deformation quantizations $(W,\F[[\hbar]],\star^1_\hbar)$ and $(W,\F[[\hbar]],\star^2_\hbar)$
of the same Poisson space $(W,\F,\{\cdot,\cdot\})$ are said to be equivalent if
there exists a $\K[[\hbar]]$-linear map $T\colon\F[[\hbar]]\rightarrow\F[[\hbar]]$ expressed
as a formal series $T=\mathrm{Id}_\mathcal{F} + \sum_{r\geq1} \hbar^r T_r$ satisfying:
\vskip2mm
\noindent i) $T_r\colon \F\rightarrow\F $, $r\geq1$, are differential operators in the sense of Def.~\ref{def:op},
vanishing on constants,
\vskip2mm
\noindent ii) $T(F)\star^1_\hbar T(G) = T(F \star^2_\hbar G),\ \forall F,G\in \F$.
\end{definition}

\section{Moyal product on a Hilbert space}\label{sec:toy}

We present an infinite-dimensional version of the Moyal product defined on a class
of smooth functions specified by a Hilbert-Schmidt type of conditions on their derivatives.

\subsection{Poisson structure}\label{sub:poisson}

Let $\H$ be an infinite-dimensional separable Hilbert space.
We consider the phase-space $W=\H\oplus\H^*$ endowed with its canonical strong symplectic
structure $\omega((x_1,\eta_1),(x_2,\eta_2))= \eta_1(x_2) - \eta_2(x_1)$, where
$x_1,x_2\in\H$ and $\eta_1,\eta_2\in\H^*$.

Let $F\colon W\rightarrow \mathbb{C}$ be a $C^\infty$ function (in the Fr\'echet sense). 
We shall denote by $D_1F(x,\eta)$ (resp. $D_2F(x,\eta)$) the first (resp. second) partial Fr\'echet derivative of $F$ 
evaluated at point $(x,\eta)\in W$. With the identification 
$\H^{**}\sim\H$ we have  $D_1F(x,\eta)\in \H^*$
and $D_2F(x,\eta)\in \H$. Let $\langle\cdot, \cdot \rangle\colon\H^* \times\H\rightarrow\K$ be the canonical
pairing between $\H$ and $\H^*$.

With these notations, the bracket associated with the 
canonical symplectic structure on $W$ takes the form:
\begin{equation}\label{poisson}
\{F,G\}(x,\eta) 
=\langle D_1F(x,\eta),D_2G(x,\eta)\rangle-\langle  D_1G(x,\eta),D_2F(x,\eta)\rangle, 
\end{equation}
where $F,G\in C^\infty(W,\K)$.
\vskip2mm
\begin{proposition}\label{prop:poisson}
The space $W$ endowed with the bracket (\ref{poisson}) is an infinite-dimen\-sional 
Poisson space or, equivalently, $(C^\infty(W,\K), \{\cdot,\cdot\})$ is a Poisson algebra.
\end{proposition}

\begin{proof} One has only to check that the map $(x,\eta)\mapsto \{F,G\}(x,\eta)$ belongs
to $C^\infty(W,\K)$ for any $F,G\in C^\infty(W,\K)$. Then Leibniz property and Jacobi identity will follow.
For $F,G\in C^\infty(W,\K)$, the maps $(x,\eta)\mapsto (D_1F(x,\eta),D_2G(x,\eta))$ and 
$(\xi, y) \mapsto \langle \xi, y \rangle$ belong to $C^\infty(W,\H^*\times\H)$
and $C^\infty(\H^*\times\H,\K)$,  respectively. The map $(x,\eta)\mapsto \{F,G\}(x,\eta)$, as composition of 
$C^\infty$ maps, is therefore in $C^\infty(W,\K)$.\hfill$\square$
\end{proof}
\vskip2mm
For any orthonormal basis $\{e_i\}_{i\geq1}$ in $\H$ and dual basis $\{e^*_i\}_{i\geq1}$ in $\H^*$,
the complex number $\partial_iF(x,\eta)$ shall denote the partial derivative of $F$ evaluated at $(x,\eta)$ 
in the direction
of $e_i$, i.e. $\partial_iF(x,\eta)=DF(x,\eta).(e_i,0)=D_1F(x,\eta).e_i$, and, similarly, 
$\partial_{i^*}F(x,\eta)=DF(x,\eta).(0,e_i^*)=D_2F(x,\eta).e_i^*$ 
is the partial derivative in the direction of $e_i^*$. Notice that $i^*$ should not be considered
as a different index from~$i$ when sums are involved, it is merely a mnemonic notation to distinguish
partial derivatives in $\H$ and in $\H^*$.

For $F\in C^\infty(W,\K)$,
we have for any $(x,\eta)\in W$ that $\sum_{i\geq1} |\partial_iF(x,\eta)|^2<\infty$ and
$\sum_{i\geq1} |\partial_{i^*}F(x,\eta)|^2<\infty$, hence
the Poisson bracket~(\ref{poisson}) admits an equivalent form in terms of
an absolutely convergent series: 
\begin{equation}\label{poissonbis}
\{F,G\}(x,\eta)  = \sum_{i\geq1} 
\big(\partial_iF(x,\eta) \partial_{{i}^*}G(x,\eta)-\partial_iG(x,\eta) \partial_{{i}^*}F(x,\eta)\big).
\end{equation}

\subsection{Functions of Hilbert-Schmidt type}\label{sub:hs}

We now define a subalgebra of $C^\infty(W,\K)$ suited for our discussion. 
Let us start with some definitions and notations.

For any $F\in C^\infty(W,\K)$ and  $(x,\eta)\in W$, the higher derivatives 
$$
D^{(r)}F(x,\eta)\colon W\times\cdots\times W\rightarrow \K,\quad r\geq 1,
$$
are bounded symmetric $r$-linear maps and  
partial derivatives of $F$ will be denoted 
$D^{(r)}_{\alpha_1 \cdots \alpha_r}F(x,\eta)$ where $\alpha_1,\ldots,\alpha_r$ are taking values $1$ or $2$.
Let us introduce:
\begin{equation}\label{defH}
\H^{(\alpha)} = \begin{cases}
\H,  & \mbox{if}\ \alpha=1;\\
\H^*  ,  & \mbox{if}\ \alpha=2.
\end{cases}
\qquad
\alpha^\flat = \begin{cases}
2 , & \mbox{if}\ \alpha=1;\\
1 ,   & \mbox{if}\ \alpha=2.
\end{cases}
\qquad
i^{(\alpha)} = \begin{cases}
i , & \mbox{if}\ \alpha=1;\\
i^* ,   & \mbox{if}\ \alpha=2.
\end{cases}
\end{equation}
Also $i^\sharp$ will stand for either $i$ or $i^*$.
With these notations, partial derivatives of $F$ are bounded $r$-linear maps: 
$$
D^{(r)}_{\alpha_1 \cdots \alpha_r}F(x,\eta)\colon \H^{(\alpha_1)}\times\cdots\times \H^{(\alpha_r)}\rightarrow \K.
$$

It is convenient to introduce new symbols such as $\partial_{ij^*k}$
for higher partial derivatives, e.g., $\partial_{ij^*k}F(x,\eta)\in \K$ stands for 
$D^{(3)}F(x,\eta).((e_i,0),(0,e_j^*),(e_k,0))=D^{(3)}_{121}F(x,\eta).(e_i,e_j^*,e_k)$, 
where $\{e_i\}_{i\geq1}$ (resp. $\{e_i^*\}_{i\geq1}$)
is an orthonormal basis in $\H$ (resp. $\H^*$). 
\vskip2mm
\begin{definition}\label{defhs} Let $\{e_i\}_{i\geq1}$ 
be an orthonormal basis in $\H$ and $\{e_i^*\}_{i\geq1}$ be the dual basis in $\H^*$. 
Functions of Hilbert-Schmidt type are functions $F$ in $C^\infty(W,\K)$ such that

\begin{equation}\label{hs}
\sum_{i_1, \ldots, i_r \geq1} |\partial_{i^{\sharp}_1\cdots i^{\sharp}_r}F(x,\eta)|^2 < \infty, 
\quad \forall r\geq1, \forall (x,\eta)\in W.
\end{equation}

The sums involved have to be  interpreted in the sense of  summable families.
By Schwarz lemma for partial derivatives, it should be understood that Eq.~(\ref{hs}) 
represents $r+1$ distinct sums corresponding to all of the choices $i^\sharp=i$ or $i^*$.
The set of functions of Hilbert-Schmidt type on $W$ will be denoted by $\F_{HS}$.
\end{definition}
\vskip2mm
The definition above is independent of the choice of the orthonormal basis.
\vskip2mm
\begin{remark} 
Let $\N_*$ be the set of positive integers. For each $r\geq1$,
the set of families of elements $\{x_I\}_{I\in \N_*^r}$ in $\K$ 
such that $\sum_{I\in \N_*^r}{|x_I|^2}<\infty$ is the
Hilbert space $\ell^2({\N_*^r)}$ for the usual operations and inner product.
Then condition~(\ref{hs}) can be equivalently stated in the following way:
$F\in C^\infty(W,\K)$ is of Hilbert-Schmidt type if and only if,
for any $r\geq1$ and any $(x,\eta)\in W$, the $2^r$ families 
$\{\partial_{i^{\sharp}_1\cdots i^{\sharp}_r}F(x,\eta)\}_{(i^{\sharp}_1,\ldots, i^{\sharp}_r)\in  \N_*^r}$
belong to $\ell^2({\N_*^r)}$. 
\end{remark}

\begin{remark}
The set $\F_{HS}$ does not contain all of the (continuous) polynomials on $W$. For
example, the polynomial $P(y,\xi)=\langle\xi, y\rangle$ is not in $\F_{HS}$
as $\sum_{i,j\geq1} |\partial_{ij^*}P(x,\eta)|^2$ $= \sum_{i,j\geq1} \delta_{ij}=\infty$.
In a quantum field theory context, the polynomial~$P$ corresponds to a free Hamiltonian
in the holomorphic representation.
\end{remark}
\vskip2mm
\begin{proposition}\label{pro:subhs}
The set of functions of Hilbert-Schmidt type is an Abelian subalgebra of $C^\infty(W,\K)$
for the pointwise product of functions.
\end{proposition}
\begin{proof}
Let $a,b\in \K$ and $F,G\in \F_{HS}$. It is clear from Remark~1 that  
$aF + bG$ is in $F,G\in \F_{HS}$. The product $FG$ belongs to
$\F_{HS}$ as a consequence of the Leibniz rule for the  derivatives and from:  
if $\{x_I\}_{I\in \N_*^r}\in \ell^2({\N_*^r)}$ and $\{y_J\}_{J\in \N_*^s}\in \ell^2({\N_*^s)}$,
then $\{x_Iy_J\}_{I\times J\in \N_*^{r+s}}\in \ell^2({\N_*^{r+s})}$. \hfill$\square$
\end{proof}
\vskip2mm
Moreover the Poisson bracket~(\ref{poisson}) restricts to $\F_{HS}$ and we have:
\vskip2mm
\begin{proposition}\label{pro:poihs}
$(W,\F_{HS}, \{\cdot,\cdot\})$ is a Poisson space.
\end{proposition}
\begin{proof}
Let $F$ and $G$ be in $\F_{HS}$. According to the proof of Prop.~\ref{prop:poisson}, the 
map $\Phi\colon (x,\eta)\mapsto \langle D_1F(x,\eta),D_2G(x,\eta)\rangle$ is in $ C^\infty(W,\K)$
and splits as follows:

\begin{equation*}
\begin{CD}
\Phi\colon W @>{\Psi_1}>>\H^*\times \H  @>{\Psi_2}>>\K
\end{CD}
\end{equation*}

$$
(x,\eta)\mapsto (D_1F(x,\eta),D_2G(x,\eta)) \mapsto \langle D_1F(x,\eta),D_2G(x,\eta)\rangle,
$$
\vskip2mm
\noindent where both  $\Psi_1$ and  $\Psi_2$ are $C^\infty$ maps. We only need to check
that $\Phi$ is of Hilbert-Schmidt type.

By applying the chain rule
to $\Phi=\Psi_2\circ \Psi_1$, it is easy to see that we can freely interchange partial derivatives
with the sum sign and we get that the partial derivatives of $\Phi$ is a 
finite sum of terms of the form:
\begin{equation}\label{summ}
a_{j_1^\sharp\cdots j_r^\sharp,k_1^\sharp\cdots k_s^\sharp}
\equiv\sum_{i\geq1} \partial_{ij_1^\sharp\cdots j_r^\sharp}F(x,\eta)  
\partial_{i^*k_1^\sharp\cdots k_s^\sharp}G(x,\eta).
\end{equation}
The Cauchy-Schwarz inequality implies that the family 
$\{ a_{j_1^\sharp\cdots j_r^\sharp,k_1^\sharp\cdots k_s^\sharp}\}$ is in $\ell^2(\N^{r+s}_*)$ and thus
$\Phi$ belongs to $\F_{HS}$. Hence $\F_{HS}$ is closed under the Poisson bracket.\hfill$\square$
\end{proof}

\subsection{Moyal star-product on $W$}\label{sub:moyal}

We are now in position to define the Moyal star-product on $W$ as an associative product
on $\F_{HS}[[\hbar]]$.

For $F,G\in \F_{HS}$, $(x,\eta)\in W$, $r\geq1$, $\alpha_1,\ldots,\alpha_r,\beta_1\ldots,\beta_r$ equal to
$1$ or $2$, and with the notations introduced previously, let us define:
\begin{equation}\label{conr}
\langle\langle D^{(r)}_{\alpha_1\cdots \alpha_r}F,D^{(r)}_{\beta_1\cdots \beta_r}G\rangle\rangle(x,\eta)
= \sum_{i_1,\ldots,i_r\geq1}
\partial_{i^{(\alpha_1)}_1\cdots i^{(\alpha_r)}_r}F(x,\eta)\
\partial_{i^{(\beta_1)}_1\cdots i^{(\beta_r)}_r}G(x,\eta).
\end{equation}
\vskip2mm
\begin{remark}
The preceding definition does not depend on the choice of the orthonormal basis in $\H$ and
the series is absolutely convergent as a consequence of the Cauchy-Schwarz inequality.
\end{remark}
\vskip2mm
Let $\Lambda$ be the canonical symplectic $2\times2$-matrix 
with $\Lambda^{12}=+1$. 
As in the finite-dimensional case, the powers of the Poisson bracket~(\ref{poisson}) are defined as:
\begin{equation}\label{eq:cr}
C_r(F,G)
 = \sum_{\alpha_1,\ldots,\alpha_r=1,2}\sum_{\beta_1,\ldots,\beta_r=1,2}
\Lambda^{\alpha_1\beta_1}\cdots \Lambda^{\alpha_r\beta_r} \langle \langle D^{(r)}_{\alpha_1\cdots \alpha_r}F,
D^{(r)}_{\beta_1\cdots \beta_r}G\rangle\rangle.
\end{equation}

The next Proposition shows that the $C_r$ are bidifferential operators in the sense of Def.~\ref{def:op}
and they close on $\F_{HS}$. We shall use a specific version of
the Hilbert tensor product $\otimes$ between Hilbert spaces 
(see e.g. Sect.~2.6 in \cite{kadison.ringrose:1997a}). Let $\H_1,\ldots,\H_r$ be Hilbert
spaces with orthonormal bases $\{e^{(1)}_i\}_{i\geq1},\ldots,\{e^{(r)}_i\}_{i\geq1}$. 
There exists a Hilbert space $\T=\H_1\otimes \cdots \otimes\H_r$ and
a bounded $r$-linear map $\Psi\colon(x_1,\ldots,x_r)\mapsto x_1\otimes \cdots \otimes x_r$
from  $\H_1\times \cdots \times\H_r$ to $\H_1\otimes \cdots \otimes\H_r$ satisfying:
$$
\sum_{i_1,\ldots,i_r\geq1} |\langle\Psi(e^{(1)}_{i_1},\ldots,e^{(r)}_{i_r}),x\rangle|^2 <\infty, \quad \forall x\in \T,
$$
such that for any bounded $r$-linear form: $ \Xi\colon \H_1\times \cdots \times\H_r\rightarrow \K$
satisfying:
$$
\sum_{i_1,\ldots,i_r\geq1} |\Xi(e^{(1)}_{i_1},\ldots,e^{(r)}_{i_r})|^2 <\infty,
$$
there exists a unique bounded linear form $L$ on $\T$ so that $\Xi=L\circ\Psi$.
This universal property allows to identify $\Xi$ to an element of $\T^*$.
\vskip2mm
\begin{proposition}\label{pro:cr}
For $F,G\in \F_{HS}$ and $r\geq1$, the map $(x,\eta)\mapsto C_r(F,G)(x,\eta)$ belongs to the
space of functions of Hilbert-Schmidt type $\F_{HS}$.
\end{proposition}

\begin{proof}
Each term in the finite sum~(\ref{eq:cr}) is of the form
\begin{equation}\label{eq:pm}
 \langle\langle D^{(r)}_{\alpha_1\cdots \alpha_r}F(x,\eta),
D^{(r)}_{\alpha^\flat_1\cdots \alpha^\flat_r}G(x,\eta)\rangle\rangle,
\end{equation}
where $\alpha_1,\ldots, \alpha_r= 1$ or $2$. From the definition of $\F_{HS}$, expression~(\ref{eq:pm})
is well defined for any $F,G\in \F_{HS}$ and thus defines a function  on $W$:
$$
\Phi\colon(x,\eta)\mapsto \langle\langle D^{(r)}_{\alpha_1\cdots \alpha_r}F(x,\eta),
D^{(r)}_{\alpha^\flat_1\cdots \alpha^\flat_r}G(x,\eta)\rangle\rangle.
$$

The case $r=1$ has been already proved in Prop.~\ref{pro:poihs}. For $r\geq2$, we need
to slightly modify the argument used in the proof of Prop.~\ref{pro:poihs}
since the bilinear map $\langle\langle\ ,\ \rangle\rangle$ defined by~(\ref{conr}) is not
a bounded bilinear form on the product of Banach spaces:
\begin{equation}\label{pairun}
\mathcal{L}^r(\H^{(\alpha_1)},\ldots,\H^{(\alpha_r)};\K)\times 
\mathcal{L}^r(\H^{(\alpha^\flat_1)},\ldots, \H^{(\alpha^\flat_r)};\K).
\end{equation}

In order to show that $\Phi$ is in $C^\infty(W,\K)$, we shall use the universal property
of  the Hilbert tensor product $\otimes$ mentioned above. For $F,G\in \F_{HS}$,
we can consider the bounded $r$-linear maps
 $ D^{(r)}_{\alpha_1\cdots \alpha_r}F(x,\eta)$ and $D^{(r)}_{\alpha^\flat_1\cdots \alpha^\flat_r}G(x,\eta)$
as bounded linear forms on  $\H^{(\alpha)}\equiv\H^{(\alpha_1)}\otimes\cdots\otimes \H^{(\alpha_r)}$
and $\H^{(\alpha^\flat)}\equiv\H^{(\alpha^\flat_1)}\otimes\cdots\otimes \H^{(\alpha^\flat_r)}$, respectively.
Then the unbounded bilinear map $\langle\langle\ ,\ \rangle\rangle$ on the product of spaces~(\ref{pairun}) 
restricts to the natural pairing of $\H^{(\alpha^\flat)}\sim \H^{(\alpha)*}$ and $\H^{(\alpha)}$
which is a smooth map. This shows that $\Psi$ belongs to $C^\infty(W,\K)$.

An argument similar to the one used in the proof of Prop.~3 shows that the partial derivatives 
of $\Phi$ involve a finite sum of terms of the form:
$$
\sum_{i_1,\ldots i_r\geq1} \partial_{i^{(\alpha_1)}_1\cdots i^{(\alpha_r)}_r j_1^\sharp\cdots j_a^\sharp}F(x,\eta)  
\ \partial_{i^{(\alpha^\flat_1)}_1\cdots i^{(\alpha^\flat_r)}_r k_1^\sharp\cdots k_b^\sharp}G(x,\eta),
$$
where we have used the notations introduced at the beginning of Subsection~3.2.
A~direct application of the Cauchy-Schwarz inequality gives:
\begin{align*}
&\sum_{j_1,\ldots, j_a\geq1}
\sum_{k_1,\ldots, k_b\geq1}
\Big|
\sum_{i_1,\ldots i_r\geq1} \partial_{i^{(\alpha_1)}_1\cdots i^{(\alpha_r)}_r j_1^\sharp\cdots j_a^\sharp}F 
\ \partial_{i^{(\alpha^\flat_1)}_1\cdots i^{(\alpha^\flat_r)}_r k_1^\sharp\cdots k_b^\sharp}G
\Big|^2 \leq\\
& \sum_{i_1,\ldots i_r,j_1,\ldots, j_a\geq1}
\big|\partial_{i^{(\alpha_1)}_1\cdots i^{(\alpha_r)}_r j_1^\sharp\cdots j_a^\sharp}F\big|^2\
 \sum_{i_1,\ldots i_r,k_1,\ldots, k_b\geq1}
\big| \partial_{i^{(\alpha^\flat_1)}_1\cdots i^{(\alpha^\flat_r)}_r k_1^\sharp\cdots k_b^\sharp}G\big|^2<\infty.
\end{align*}
This shows that  $\Phi$ is of Hilbert-Schmidt type and hence $(x,\eta)\mapsto C_r(F,G)(x,\eta)$
belongs to $\F_{HS}$.\hfill$\square$
\end{proof}
\vskip2mm
We summarize all the previous facts in the following:
\vskip2mm
\begin{theorem}
Let the $C_r$'s be given by (\ref{eq:cr}), then the formula
\begin{equation}
F\star^M_\hbar G = F G + \sum_{r\geq1}\frac{\hbar^r}{r!} C_r(F,G), 
\end{equation} 
defines an associative product on $\F_{HS}[[\hbar]]$ and hence a
 deformation quantization $(W,\F_{HS}[[\hbar]],\star^M_\hbar)$ 
of the Poisson space $(W,\F_{HS},\{\cdot,\cdot\})$.
\end{theorem}

\begin{proof}
That $\star^M_\hbar\colon\F_{HS}[[\hbar]]\times \F_{HS}[[\hbar]]\rightarrow \F_{HS}[[\hbar]]$
is a bilinear map is a direct consequence of Prop.~\ref{pro:cr} and Prop.~\ref{pro:subhs}. 
Associativity follows from the same combinatorics  used in the finite-dimensional case and the
fact that the derivatives distribute in the pairing $\langle\langle\ ,\ \rangle\rangle$  
defining the $C_r$'s.\hfill$\square$
\end{proof}

\section{On the equivalence of deformation quantizations on $W$}

We end this article by a discussion on the issue of equivalence of star-products on $W$.
In the finite-dimensional case (i.e. on $\R^{2n}$ endowed with its canonical Poisson bracket),
it is well known that all star-products are equivalent to each other. The situation we are dealing
with here, although it is a direct generalization of the flat finite-dimensional case, allows
inequivalent deformation quantizations. We will illustrate this fact on a family
of star-products of exponential type containing the important case of the normal star-product.

\subsection{The Hochschild complex}

The space of functions of Hilbert-Schmidt type $\F_{HS}$  being an associative
algebra over $\K$, we can consider the Hochschild complex 
$C^\bullet(\F_{HS},\F_{HS})$ and its cohomology $H^\bullet(\F_{HS},\F_{HS})$.

One has first to specify a class of cochains that would define the Hochschild complex.
Here the cochains are simply  $r$-differential operators in the sense of Def.~\ref{def:op}
which vanishes on constants. The case where $\F=\F_{HS}$ in Def.~\ref{def:op}
allows a more precise description of the $r$-differential operators. Consider
an $r$-differential operator defined by:
$$
\mathbf{A}(F_1,\ldots,F_r)(w)=\sum_{n_1,\ldots,n_r\geq 0} 
a^{(n_1,\cdots ,n_r)}(w).(D^{(n_1)}F_1(w),\cdots,D^{(n_r)}F_r(w)),
$$
where $F_1,\ldots,F_r\in \F_{HS}$ and $w=(x,\eta)\in W$.
For $F\in \F_{HS}$, the higher derivative $D^{(m)}F(x,\eta)$ defines an element
of the $m^{\rm th}$ tensor power of $W$ (here we identify $W^*$ with $W$). 
For a fixed $w=(x,\eta)\in W$ and $m\geq0$, the linear map $F\mapsto D^{(m)}F(x,\eta)$ from 
$\F_{HS}$ to $W^{\stackrel{m}{\otimes}}$ is onto and we can look at the restriction of the
$r$-linear form
$$
a^{(n_1,\ldots,n_r)}(w)\colon \mathcal{D}_w^{(n_1,\ldots,n_r)}\subset
\mathcal{L}^{n_1}_{\rm sym}(W,\K)\times\cdots\times \mathcal{L}^{n_r}_{\rm sym}(W,\K)\rightarrow \K,
$$
to the product $W^{\stackrel{n_1}{\otimes}}\times\cdots\times W^{\stackrel{n_r}{\otimes}}$ 
as a bounded $r$-linear form:
$$
\tilde a^{(n_1,\ldots,n_r)}(w)\colon W^{\stackrel{n_1}{\otimes}}\times\cdots\times W^{\stackrel{n_r}{\otimes}}
\rightarrow \K,
$$
such that $w\mapsto \tilde a^{(n_1,\ldots,n_r)}(w)$ is a smooth map from $W$ to  
$\mathcal{L}^{r}(W^{\stackrel{n_1}{\otimes}}\times\cdots\times W^{\stackrel{n_r}{\otimes}},\K)$.

The Leibniz rule for the derivatives of a product can be written here for $F,G\in  \F_{HS}$ as:
\begin{align*}
D^{(m)}&(FG)(w).(w_1,\ldots,w_m) \cr
&= \frac{1}{m!} \sum_{\sigma\in \mathfrak{S}_m} \sum_{k=0}^{k=m} 
\binom{m}{k}\ (D^{(k)}F(w) \otimes D^{(m-k)}G(w)) (w_{\sigma_1}\otimes\cdots\otimes w_{\sigma_m}),
\end{align*}
where $\mathfrak{S}_m$ is the symmetric group of degree~$m$. 

Let 
$\mathbf{A}(F)=\sum_{m\geq0} \tilde a^{(m)}(w).(D^{(m)}F(w))$
be a differential operator with 
$\tilde a^{(m)}\in C^\infty(W,\mathcal{L}(W^{\stackrel{m}{\otimes}},\K))$, 
then it follows from the above form for the Leibniz rule that for $F,G\in \F_{HS}$,
$\mathbf{A}(FG)$ can be written as a finite sum
$$
\sum_{m_1,m_2\geq0} \tilde b^{(m_1,m_2)}(w).(D^{(m_1)}F(w),D^{(m_2)}G(w))
$$ for some
$\tilde b^{(m_1,m_2)}\in C^\infty(W,\mathcal{L}^2(W^{\stackrel{m_1}{\otimes}}\times W^{\stackrel{m_2}{\otimes}},\K))$,
and thus $(F,G)\mapsto \mathbf{A}(FG)$ is a bidifferential operator. The generalization to $r$-differential operators
of this fact is straightforward.

The Hochschild complex consists of multidifferential operators 
vanishing on constants,
i.e., $C^\bullet(\F_{HS},\F_{HS})=\oplus_{k\geq0} C^k(\F_{HS},\F_{HS})$ 
where, for $k\geq1$, the space of \mbox{$k$-cochains} is: 
\begin{align*}
&C^k(\F_{HS},\F_{HS})\\
&=\{\mathbf{A}\colon \F_{HS}^{k}\rightarrow \F_{HS}\ | \
\mathbf{A}\ {\rm is\ a}\ k{\rm -differential\ operator\ vanishing\ on}\ \K\}.
\end{align*}

The differential of a $k$-cochain $\mathbf{A}$ is the  $(k+1)$-linear map  $\delta\mathbf{A}$ given by:
\begin{eqnarray} \label{bH}
\delta\mathbf{A}(F_0,\ldots,F_k)&=&F_0\mathbf{A}(F_1,\ldots,F_k)-\mathbf{A}(F_0F_1,F_2,\ldots,F_k)+
\cdots\\
&+&(-1)^k\mathbf{A}(F_0,F_1,\ldots,F_{k-1}F_k)+(-1)^{k+1}\mathbf{A}(F_0,\ldots,
F_{k-1})F_k.\nonumber
\end{eqnarray}
satisfies $\delta^2=0$, and according to the discussion above 
$\delta\mathbf{A}$ is a $(k+1)$-differential operator,
vanishing on constants whenever $\mathbf{A}$ does. Thus $\delta\mathbf{A}$
is a cochain that belongs to $C^{k+1}(\F_{HS},\F_{HS})$, 
hence we indeed have a complex.

A $k$-cochain $\mathbf{A}$ is a $k$-cocycle if $\delta\mathbf{A}=0$. 
We denote by $Z^k(\F_{HS},\F_{HS})$ the space of $k$-cocycles and by $B^k(\F_{HS},\F_{HS})$ the
space of those $k$-cocycles which are coboundaries.
The $k^{\rm th}$ Hochschild cohomology space of $\F_{HS}$ valued in $\F_{HS}$
is defined as the quotient $H^k(\F_{HS},\F_{HS})=Z^k(\F_{HS},\F_{HS})/B^k(\F_{HS},\F_{HS})$.

\subsection{Star-products of the exponential type}

Let $\mathcal{B}(\H)$ denote the algebra of bounded operators on $\H$ and $\mathcal{B}_2(\H)$, the
two-sided $*$-ideal of Hilbert-Schmidt operators on $\H$ . 
We shall describe a family of deformation quantizations  
$\{(W,\F_{HS}[[\hbar]],\star^A_\hbar)\}_{A\in \mathcal{B}(H)}$.
Each star-product~$\star^A_\hbar$ where $A\in \mathcal{B}(H)$ shall 
be the exponential of a Hochschild $2$-cocycle, with the Moyal star-product
corresponding to the case $A=0$. It will turn out that the set of
 equivalence classes of star-products of this type is parameterized by 
$\mathcal{B}(\H)/\mathcal{B}_2(\H)$.

Let $A \in \mathcal{B}(\H)$. For $F,G\in \F_{HS}$, the map 
\begin{equation}\label{eadef}
(x,\eta)\mapsto E_A(F,G)(x,\eta)\equiv\langle D_1F(x,\eta), A D_2G(x,\eta)\rangle
+ \langle D_1G(x,\eta), A D_2F(x,\eta)\rangle
\end{equation}
defines a smooth function on $W$ and is symmetric in $F$ and $G$. Moreover we have:
\begin{proposition}\label{lemma:ea}
The bilinear map $(F,G)\mapsto  E_A(F,G)$ is  a Hochschild $2$-cocycle.
\end{proposition}

\begin{proof}
It is clear that $\delta E_A=0$ and, since $E_A$ vanishes on constants, it is sufficient to check 
that $E_A$ is a bidifferential operator on $(W,\F_{HS})$, i.e., that the smooth function 
$(x,\eta)\mapsto E_A(F,G)(x,\eta)$ is of Hilbert-Schmidt type for any $F,G\in \F_{HS}$.
Let $\{e_i\}_{i\geq 1}$ be an orthonormal basis of $\H$ and $\{e^*_i\}_{i\geq 1}$ the dual basis.
Consider the basis $\{f_i\}_{i\geq 1}$ in $W=\mathcal{H}\oplus \mathcal{H}^*$ defined by
$f_i=(e_{\frac{i+1}{2}},0)$ if $i$ is odd, and $f_i=(0,e^*_{\frac{i}{2}})$ if $i$ is even.
To show that  $E_A(F,G)$ is in $\F_{HS}$ is equivalent to show that  
\begin{equation}\label{eahs}
\sum_{i_1,\ldots,i_n\geq1}| D^{(n)}E_A(F,G)(x,\eta).(f_{i_1},\ldots,f_{i_n})|^2 <\infty.
\end{equation}
holds for any $n\geq1$ and $w=(x,\eta)\in W$.

The chain rule applied to
\begin{equation*}
\begin{CD}
\Phi\colon W @>{\Psi_1}>>\H^*\times \H  @>{\Psi_2}>>\K
\end{CD}
\end{equation*}

$$
(x,\eta)\mapsto (D_1F(x,\eta),D_2G(x,\eta)) \mapsto \langle D_1F(x,\eta),AD_2G(x,\eta)\rangle,
$$
\vskip2mm
\noindent shows that the derivatives of $E_A$ distributes in the pairing $\langle\ , \ \rangle$.

The $n^{\rm th }$ derivative of $E_A(F,G)$ in the direction of $(f_{i_1},\ldots,f_{i_n})$ is
a finite sum of terms of the form:
\begin{equation}\label{da}
\langle D^{(r)}D_1F(x,\eta).(f_{k_1},\ldots,f_{k_r}), A D^{(s)}D_2G(x,\eta).(f_{l_1},\ldots,f_{l_s})\rangle,
\end{equation}
where $r+s=n$ and $(k_1,\ldots,k_r, l_1,\ldots,l_s)$ is a permutation of $(i_1,\ldots,i_n)$, and
similar terms  with $F$ and $G$ inverted.
It is worth noting that $D^{(r)}D_1F(x,\eta).(f_{k_1},\ldots,f_{k_r})$ is the element of $\H^*$ defined
by the bounded linear form 
$$h\mapsto D^{(r+1)}F(x,\eta).((h,0),f_{k_1},\ldots,f_{k_r})$$
 on $\H$
and, similarly, 
$D^{(s)}D_2G(x,\eta).(f_{l_1},\ldots,f_{l_s})$ is the element of $\H$ defined by the bounded
linear form on $\H^*$: 
$$
\xi\mapsto D^{(s+1)}G(x,\eta).((0,\xi),f_{l_1},\ldots,f_{l_s}).$$
The modulus squared of (\ref{da}) is bounded by the constant
$
||A||^2 \alpha_{k_1\ldots k_r} \beta_{l_1\ldots l_s}
$, where
\begin{align*}
\alpha_{k_1\ldots k_r}=&\sum_{i\geq1}|D^{(r+1)}F(x,\eta).((e_i,0),f_{k_1},\ldots,f_{k_r})|^2,\\
\beta_{l_1\ldots l_s}=&\sum_{j\geq1}|D^{(s+1)}G(x,\eta).((0,e^*_j),f_{l_1},\ldots,f_{l_s})|^2.
\end{align*}
Since $F,G\in\F_{HS}$, we have $\ds\sum_{k_1,\ldots,k_r\geq1}\alpha_{k_1\ldots k_r}<\infty$ and
$\ds\sum_{l_1,\ldots,l_s\geq1}\beta_{l_1\ldots l_s}<\infty$, from which inequality~(\ref{eahs})
follows.\hfill$\square$
\end{proof}

For any $A \in \mathcal{B}(\H)$, let us define:
\begin{align}\label{coneafg}
C^A_1(F,G) & = \{F,G\} + E_A(F,G)\\
           & = \langle D_1F,(A + \I)D_2G\rangle + \langle D_1G,(A - \I)D_2F\rangle\nonumber,
\end{align}
where $\I$ is the identity operator on $\H$. Plainly, $C^A_1$ is a $2$-cocycle with constant coefficients,
and one can define a star-product by taking the exponential of $C^A_1$:
$$
F\star^A_\hbar G = \exp(\hbar C^A_1)(F,G)= FG + \sum_{r\geq 1} \frac{\hbar^r}{r!}  C^A_r(F,G),
$$
where $C^A_r= (C^A_1)^r$ in the sense of bidifferential operators. This formula 
defines an associative product on $\F_{HS}[[\hbar]]$, and
we get a  family of deformation quantizations  
$\{(W,\F_{HS}[[\hbar]],\star^A_\hbar)\}_{A\in \mathcal{B}(H)}$ of $(W,\F_{HS},\{,\cdot,\cdot\})$.
This family of star-products is easily described by their symbols. 
Let us consider the following family of smooth functions on $W$:
\begin{equation}\label{phiyxi}
\Phi_{y,\xi}(x,\eta) =\exp(\langle\eta,y \rangle + \langle\xi,x \rangle), \quad
x,y\in \H,\ \eta,\xi\in \H^*.
\end{equation}
The $\Phi_{y,\xi}$'s belong to $\F_{HS}$ and from (\ref{coneafg}) we deduce that:
$$
C^A_r(\Phi_{y,\xi},\Phi_{y',\xi'}) = 
(\langle\xi, (A + \I)y' \rangle + \langle\xi',(A - \I)y \rangle)^r\ \Phi_{y+y',\xi+\xi'},
$$
and consequently:
\begin{equation}\label{asymbol}
\Phi_{y,\xi}\star^A_\hbar \Phi_{y',\xi'} = 
\exp\big(\hbar(\langle\xi, (A + \I)y' \rangle + \langle\xi',(A - \I)y \rangle)\big)\ \Phi_{y+y',\xi+\xi'}.
\end{equation}

\begin{example}
Set $A=\I$ in (\ref{coneafg}), then $C^\I_1(F,G) = 2 \langle D_1F,D_2G\rangle$ and 
the corresponding star-product reads 
\begin{equation}\label{normal}
F\star^\I_\hbar G= F G+\sum_{r\geq1}\frac{(2\hbar)^r}{r!} 
\langle D^{(r)}_{1\cdots 1}F,D^{(r)}_{2\cdots 2}G\rangle.
\end{equation}
It is the well-known normal star-product
(or Wick or standard depending on the interpretation of the variables $(x,\eta)\in W$).
The cochains defining the normal star-product $C^\I_r$ correspond to only one term
in the sum defining the $r^{\rm th}$ cochain of the Moyal star-product, 
namely the term corresponding
to $\alpha_1=\cdots=\alpha_r=1$ and $\beta_1=\cdots=\beta_r=2$
in the sum  (\ref{eq:cr}). One would expect that conditions~(\ref{hs}) defining 
the functions of Hilbert-Schmidt type are not all needed in order to make sense
of the normal star-product and would guess that this product can be defined
on a wider class of functions.
Actually, the normal star-product defines a deformation quantization on a larger
space of functions (containing the free Hamiltonian) that we shall describe 
in a forthcoming paper.
\end{example}

At this stage, a natural question arises: are the deformation quantizations
$\{(W,\F_{HS}[[\hbar]],\star^A_\hbar)\}_{A\in \mathcal{B}(H)}$
equivalent to each other? The answer is given in the:

\begin{proposition}\label{coboundary}
Let $A \in \mathcal{B}(\H)$. The Hochschild $2$-cocycle $E_A$ defined by (\ref{eadef}) is a coboundary
if and only if $A$ is a Hilbert-Schmidt operator.
\end{proposition}
\begin{proof}
Let $A$ be a Hilbert-Schmidt operator on $\H$. The map from $W$ to $\K$ defined by
$(x,\eta)\mapsto \langle\eta,A x \rangle$ is of Hilbert-Schmidt type and therefore
there exists a bounded linear form $\tilde A\colon \H^*\otimes\H\rightarrow \K$
so that
$\langle\eta,A x \rangle = \langle\tilde A, \eta\otimes x \rangle$.
The $2$-cocycle $E_A$ can then be written as:
$$
E_A(F,G)= \langle\tilde A, D_1F\otimes D_2G + D_1G\otimes D_2F \rangle,\quad F,G\in \F_{HS}.
$$
For $F \in \F_{HS}$, the mixed derivative $D^{(2)}_{12}F$ belongs to 
$C^\infty(W, \mathcal{L}(\H \otimes \H^*,\K))\sim C^\infty(W, \H^* \otimes \H)$ and
$T_A(F)(x,\eta) = - \langle\tilde A, D^{(2)}_{12}F(x,\eta) \rangle$ defines a differential
operator on $(W, \F_{HS})$ vanishing on constants. For any two functions in $\F_{HS}$,
the Leibniz rule reads:
$$
D^{(2)}_{12}(FG)
= F D^{(2)}_{12}(G) + G D^{(2)}_{12}(F) + D_1F\otimes D_2G + D_1G\otimes D_2F,
$$
it belongs to $C^\infty(W, \H^* \otimes \H)$ 
and a simple computation gives $\delta T_A= E_A$, therefore $E_A$ is a coboundary if $A$
is in the Hilbert-Schmidt class.

Conversely, if $E_A$ is a coboundary, there exists a differential operator $S$
on $(W, \F_{HS})$ vanishing on constants, so that $E_A=\delta S$. If a term
of degree one occurs in $S$ (a derivation) it can be subtracted without 
changing $\delta S$, hence we can assume that $S$ has the form:
$$
S(F)(x,\eta)= \sum_{m\geq2} a^{(m)}(x,\eta).(D^{(m)}F(x,\eta)),
$$
where $a^{(m)}\in C^\infty(W, \mathcal{L}(W^{\stackrel{m}{\otimes}},\K))$ and
only finitely many of them are nonzero. By computing $\delta S$ and using 
$E_A=\delta S$, we find that only the term of degree $2$ contributes:
\begin{align*}
\langle & D_1F(w)\,A D_2G(w)\rangle +\langle D_1G(w),A D_2F(w) \rangle\\
&=- a^{(2)}(w).(DF(w)\otimes DG(w) + DG(w)\otimes DF(w)),\quad w=(x,\eta)\in W.
\end{align*}
If we evaluate the equality above on $F(x,\eta)=\langle \xi,x\rangle$, $\xi\in \H^*$, 
and $G(x,\eta)=\langle \eta,y\rangle$, $y\in \H$, we find (with a slight abuse of notations):
$$
\langle \xi,A y\rangle = - a^{(2)}(w).(\xi\otimes y),\quad \forall \xi\in\H^*, \forall y\in\H,
$$
from which follows that $A\in \mathcal{B}(\H)$ is a Hilbert-Schmidt operator on $\H$, 
as $a^{(2)}(w)$ is a bounded linear form on $\H^*\otimes \H$.\hfill$\square$
\end{proof}

As an  immediate consequence of Prop.~(\ref{coboundary}), we deduce a classification result
for deformation quantizations of exponential type 
$\{(W,\F_{HS}[[\hbar]],\star^A_\hbar)\}_{A\in \mathcal{B}(H)}$.

\begin{theorem}
Let $A,B\in \mathcal{B}(H)$. Two deformation quantizations $(W,\F_{HS}[[\hbar]],\star^A_\hbar)$
and $(W,\F_{HS}[[\hbar]],\star^B_\hbar)$
are equivalent if and only if $A-B$ is in the Hilbert-Schmidt class. Consequently,
the set of equivalence classes of $\{(W,\F_{HS}[[\hbar]],\star^A_\hbar)\}_{A\in \mathcal{B}(H)}$
is parameterized by $\mathcal{B}(\H)/\mathcal{B}_2(\H)$.
\end{theorem}
\begin{proof}
Suppose that $\star^A_\hbar$ and $\star^B_\hbar$ are equivalent, i.e., there exists
a formal series of differential operators vanishing on constants: 
$T=\mathrm{Id}_{\F_{HS}} + \sum_{r\geq1} \hbar^r T_r$, so that $T(F\star^A_\hbar G)= TF \star^B_\hbar TG$.
Then it follows that $C_1^A=C_1^B + \delta T_1$  and, from the definitions (\ref{eadef}) and (\ref{coneafg})
of $E_A$ and $C_1^A$, we have $E_{A-B}= E_A - E_B = \delta T_1$, showing that $E_{A-B}$ is a coboundary
and hence $A-B$ is a Hilbert-Schmidt operator on $\H$.

Conversely, if $S\equiv A-B$ is a Hilbert-Schmidt operator, it defines a bounded
linear form $\tilde S$  on $\H^*\otimes \H$ and a differential operator 
$T_1(F) = - \langle\tilde S, D^{(2)}_{12}F \rangle$ on 
$(W,\F_{HS})$ (cf. the proof of Prop.~\ref{coboundary}).
Since the star-products $\star^A_\hbar$ and $\star^B_\hbar$ are defined by
constant coefficient bidifferential operators, it is sufficient to establish
the equivalence at the level of symbols.
Now define the formal series of differential operators $T=\exp(\hbar T_1)$. 
Its symbol is given by 
$T(\Phi_{y,\xi})= \exp(\hbar \langle \xi,(B-A)y \rangle)\Phi_{y,\xi}$,
where $\Phi_{y,\xi}$ has been defined in (\ref{phiyxi}).
Using the symbol~(\ref{asymbol}) associated to a star-product, we find:
$$
T(\Phi_{y,\xi}\star^A_\hbar \Phi_{y',\xi'})=T\Phi_{y,\xi} \star^B_\hbar T\Phi_{y',\xi'},
\quad y,y'\in \H, \xi,\xi'\in \H^*.
$$
Therefore the deformation quantizations $(W,\F_{HS}[[\hbar]],\star^A_\hbar)$ and
$(W,\F_{HS}[[\hbar]],\star^B_\hbar)$ are equivalent.\hfill$\square$
\end{proof}

The Moyal and normal star-products correspond to $A=0$ and  $A=\I$ in (\ref{coneafg}), respectively.
Since the identity operator on the infinite-dimensional Hilbert space $\H$ is not in the
Hilbert-Schmidt class, we have:
\begin{corollary}
The Moyal and normal star-products are not equivalent deformations on $(W,\F_{HS}[[\hbar]])$.
\end{corollary}

\begin{remark}
One can generalize the class of exponential type of star-products 
by allowing formal series with coefficients in $\mathcal{B}{(\H)}$ in (\ref{coneafg}) or
(\ref{asymbol}). The set of equivalence classes would then be
$(\mathcal{B}(\H)/\mathcal{B}_2(\H))[[\hbar]]$. Since we did not show that any star-product
on $(W,\F_{HS},\{\cdot,\cdot\})$ is equivalent to a star-product of the exponential type,
$(\mathcal{B}(\H)/\mathcal{B}_2(\H))[[\hbar]]$ can only be considered as a
 \textit{lower bound} for the classification space of all the star-products on 
$(W,\F_{HS},\{\cdot,\cdot\})$.
\end{remark}

\begin{remark}
Recall that all the star-products on $\R^{2n}$ endowed with its canonical Poisson
structure are equivalent to each other. This fact should be put in relation with Von Neumann's
uniqueness theorem on the irreducible (continuous) representations of
Weyl systems associated to the canonical  commutation relations (CCR) $\{q_i,p_j\}=\delta_{ij}$.
The inequivalent representations of Weyl systems associated to
the infinite dimensional CCR have been described long time ago
by G\aa rding and Wightman \cite{GW}, and Segal \cite{Se}. Here the existence
of inequivalent deformation quantizations conveys the idea that there should be
a close link  between the set of equivalence classes of star-products and 
representations of Weyl systems associated to the CCR. It might turn out that they
actually are identical.
\end{remark}

\section*{Acknowledgments}
The author would like to thank Y. Maeda and S. Watamura, the organizers of the wonderful
meeting  \textit{International Workshop on Quantum Field Theory 
and Noncommutative Geometry} held in Sendai (November 2002), for their warmest hospitality.

\end{document}